\documentclass[12pt,thmsa]{article}%
\usepackage{sw20bams}
\usepackage{amsmath}
\usepackage{amsfonts}
\usepackage{amssymb}
\usepackage{graphicx}%
\providecommand{\U}[1]{\protect\rule{.1in}{.1in}}
\begin{document}

\author{Steven Finch}
\title{Random Gaussian Tetrahedra}
\date{March 21, 2022}
\maketitle

\begin{abstract}
Given independent normally distributed points $A,B,C,D$ in Euclidean
$3$-space, let $Q$ denote the plane determined by $A,B,C$ and $\tilde{D}$
denote the orthogonal projection of $D$ onto $Q$. The probability that the
tetrahedron $ABCD$ is acute remains intractable. We make some small progress
in resolving this issue. Let $\Gamma$ denote the convex cone in $Q$ containing
all linear combinations $A+r\,(B-A)+s\,(C-A)$ for nonnegative $r$, $s$. We
compute the probability that $\tilde{D}$ falls in $(B+C)-\Gamma$ to be
$0.681...$, but the probability that $\tilde{D}$ falls in $\Gamma$ to be
$0.683...$. The intersection of these two cones is a parallelogram in $Q$
twice the area of the triangle $ABC$. Among other issues, we mention the
distribution of random solid angles and sums of these.

\end{abstract}

\footnotetext{Copyright \copyright \ 2010, 2015, 2022 by Steven R. Finch. All rights
reserved.}Let $a_{1}$, $a_{2}$, $a_{3}$, $b_{1}$, $b_{2}$, $b_{3}$, $c_{1}$,
$c_{2}$, $c_{3}$, $d_{1}$, $d_{2}$, $d_{3}$ be independent normally
distributed random variables with mean $0$ and variance $1$. The points
\[%
\begin{array}
[c]{ccccccc}%
A=\left(
\begin{array}
[c]{c}%
a_{1}\\
a_{2}\\
a_{3}%
\end{array}
\right)  , &  & B=\left(
\begin{array}
[c]{c}%
b_{1}\\
b_{2}\\
b_{3}%
\end{array}
\right)  , &  & C=\left(
\begin{array}
[c]{c}%
c_{1}\\
c_{2}\\
c_{3}%
\end{array}
\right)  , &  & D=\left(
\begin{array}
[c]{c}%
d_{1}\\
d_{2}\\
d_{3}%
\end{array}
\right)
\end{array}
\]
constitute the vertices of a tetrahedron in Euclidean $3$-space. The
tetrahedron $ABCD$ is \textbf{acute} if each of its six internal dihedral
angles are less than $\pi/2$. It is known \cite{ESU} that $ABCD$ is acute if
and only if the orthogonal projection of each vertex onto the plane of the
opposite face lies within that face. Neither characterization suggests an easy
approach to finding the probability that random $ABCD$ is acute. We will
examine a variation of the latter characterization, focusing on the point $D$.
Let $Q$ denote the plane determined by $A,B,C$ and $\tilde{D}$ denote the
orthogonal projection of $D$ onto $Q$. While calculating the probability that
$\tilde{D}$ falls within the triangle $ABC$ seems difficult, we succeed in
computing the probability that $\tilde{D}$ falls in either of two convex cones
in $Q$ containing $ABC$. In fact, the cones contain the parallelogram with
vertices $A,B,C,-A+B+C$; the fourth vertex is clearly the vector sum
$(B-A)+(C-A)$ displaced so that it emanates from the point $A$. Moreover, the
intersection of the cones gives precisely the parallelogram. Further research
might uncover other probabilities for related cones in $Q$, and thus the
inclusion-exclusion principle might yield the probability that $\tilde{D}$
falls within the parallelogram. This task is best left to someone else (!) but
we hope that our work provides some inspiration along the way.

In section 1, we review a well-known proof that planar triangles $ABC$ are
acute with probability $1/4$. Returning to $3$-space, we compute in section 2
the probability that $\tilde D$ falls in the cone
\[
(-A+B+C)-\left\{  r\,(B-A)+s\,(C-A):r\geq0,\,s\geq0\right\}
\]
using the Krishnaiah bivariate $F$-ratio distribution \cite{Kr, KA, Kn}. This
is the cone with vertex $-A+B+C$ and outgoing edges parallel to vectors $A-B$,
$A-C$. In the following section, we compute the probability that $\tilde D$
falls in the cone
\[
\Gamma=A+\left\{  r\,(B-A)+s\,(C-A):r\geq0,\,s\geq0\right\}
\]
using a forgotten result of Miller's \cite{Mi}, extended via convolution. This
is the cone with vertex $A$ and outgoing edges parallel to vectors $B-A$,
$C-A$. The two probabilities are numerically close but not equal (at least not
in our view). A\ rigorous proof of such an inequality is open.

\section{Triangles and Tetrahedra}

For the moment,
\[%
\begin{array}
[c]{ccccc}%
A=\left(
\begin{array}
[c]{c}%
a_{1}\\
a_{2}%
\end{array}
\right)  , &  & B=\left(
\begin{array}
[c]{c}%
b_{1}\\
b_{2}%
\end{array}
\right)  , &  & C=\left(
\begin{array}
[c]{c}%
c_{1}\\
c_{2}%
\end{array}
\right)
\end{array}
\]
are vertices of a random Gaussian triangle in Euclidean $2$-space. Let $Q$
denote the line determined by $A,B$ and $\tilde C$ be analogous to before.
Also, let $\alpha,\beta,\gamma$ denote the angles at $A,B,C$ respectively. The
orthogonal projection of $\vec w=C-A$ onto the subspace spanned by $\vec
v=B-A$ is
\[
\frac{\vec v\cdot\vec w}{\vec v\cdot\vec v}\vec v=\frac{(B-A)\cdot
(C-A)}{(B-A)\cdot(B-A)}(B-A)
\]
and clearly $\tilde C$ falls between $A$ and $B$ on $Q$ if and only if
\[
0<\frac{(B-A)\cdot(C-A)}{(B-A)\cdot(B-A)}<1\text{.}
\]
Let
\[
p=\operatorname*{P}\left(  (B-A)\cdot(C-A)>0\right)  =\operatorname*{P}\left(
\frac{(B-A)\cdot(C-A)}{(B-A)\cdot(B-A)}>0\right)  ,
\]
then
\begin{align*}
\operatorname*{P}\left(  \frac{(B-A)\cdot(C-A)}{(B-A)\cdot(B-A)}<1\right)   &
=\operatorname*{P}\left(  (B-A)\cdot(C-A)<(B-A)\cdot(B-A)\right) \\
&  =\operatorname*{P}\left(  (B-A)\cdot\left[  (C-A)-(B-A)\right]  <0\right)
\\
&  =\operatorname*{P}\left(  (B-A)\cdot(C-B)<0\right) \\
&  =\operatorname*{P}\left(  (A-B)\cdot(C-B)>0\right)  =p
\end{align*}
by symmetry, hence
\[
\operatorname*{P}\left(  0<\frac{(B-A)\cdot(C-A)}{(B-A)\cdot(B-A)}<1\right)
=p-(1-p)=2p-1.
\]
To compute $p$, note that $(B-A)\cdot(C-A)$ is equal to
\begin{align*}
&  \frac32\left(  -\sqrt{\frac23}A+\frac B{\sqrt{6}}+\frac C{\sqrt{6}}\right)
\cdot\left(  -\sqrt{\frac23}A+\frac B{\sqrt{6}}+\frac C{\sqrt{6}}\right)
-\frac12\left(  -\frac B{\sqrt{2}}+\frac C{\sqrt{2}}\right)  \cdot\left(
-\frac B{\sqrt{2}}+\frac C{\sqrt{2}}\right) \\
&  =\frac32\left\|  -\sqrt{\frac23}A+\frac B{\sqrt{6}}+\frac C{\sqrt{6}%
}\right\|  ^{2}-\frac12\left\|  -\frac B{\sqrt{2}}+\frac C{\sqrt{2}}\right\|
^{2}\\
&  =\frac32\chi_{2}^{2}-\frac12\bar\chi_{2}^{2},
\end{align*}
a linear combination of independent chi-square distributed variables (each
with $2$ degrees of freedom). Therefore
\[
p=\operatorname*{P}\left(  \frac32\chi_{2}^{2}-\frac12\bar\chi_{2}%
^{2}>0\right)  =\operatorname*{P}\left(  \frac{\chi_{2}^{2}}{\bar\chi_{2}^{2}%
}>\frac13\right)  =\operatorname*{P}\left(  F_{2,2}>\frac13\right)  =\frac34
\]
by use of the $F$-ratio distribution with $(2,2)$ degrees of freedom. It
follows that $\tilde C$ falls between $A$ and $B$ on $Q$ with probability
$2(3/4)-1=1/2$. Equivalently, $\operatorname*{P}(\alpha>\pi/2)=1/4$ because
$\cos(\alpha)<0$ precisely when $(B-A)\cdot(C-A)<0$. Since at most one angle
of a triangle can be obtuse, we deduce that
\begin{align*}
\operatorname*{P}\left(  \text{triangle }ABC\text{ is acute}\right)   &
=1-\operatorname*{P}\left(  \text{triangle }ABC\text{ is obtuse}\right) \\
&  =1-\left(  \operatorname*{P}(\alpha>\pi/2)+\operatorname*{P}(\beta
>\pi/2)+\operatorname*{P}(\gamma>\pi/2)\right) \\
&  =1-3/4=1/4
\end{align*}
as was to be shown. Our proof imitates Eisenberg \&\ Sullivan's \cite{EiS}
approach, although we avoid angles $\alpha,\beta,\gamma$ until the last step,
using $\tilde C$ instead. Portnoy's \cite{Por} argument employs triangle
medians rather than orthogonal projections.

Let us now return to random Gaussian tetrahedra in $3$-space. Our argument is
similar but more complicated. The orthogonal projection of $\vec w=D-A$ onto
the subspace spanned by $\vec u=B-A$, $\vec v=C-A$ is
\[
\frac{\vec u\cdot\vec w}{\vec u\cdot\vec u}\vec u+\frac{\vec v\cdot\vec
w}{\vec v\cdot\vec v}\vec v=\frac{(B-A)\cdot(D-A)}{(B-A)\cdot(B-A)}%
(B-A)+\frac{(C-A)\cdot(D-A)}{(C-A)\cdot(C-A)}(C-A)
\]
and clearly $\tilde D$ falls within the desired parallelogram on $Q$ if and
only if
\[%
\begin{array}
[c]{ccc}%
0<\dfrac{(B-A)\cdot(D-A)}{(B-A)\cdot(B-A)}<1 & \text{and} & 0<\dfrac
{(C-A)\cdot(D-A)}{(C-A)\cdot(C-A)}<1.
\end{array}
\]
This is the same as requiring that
\[%
\begin{array}
[c]{ccc}%
(A-B)\cdot(D-B)>0, &  & (B-A)\cdot(D-A)>0,\\
(A-C)\cdot(D-C)>0, &  & (C-A)\cdot(D-A)>0.
\end{array}
\]
Each product is of the form $(3/2)\chi_{3}^{2}-(1/2)\bar\chi_{3}^{2}$ and
jointly they give rise to corresponding $F_{3,3}$ ratios:
\[%
\begin{array}
[c]{ccc}%
\dfrac{\left\|  \dfrac A{\sqrt{6}}-\sqrt{\dfrac23}B+\dfrac D{\sqrt{6}%
}\right\|  ^{2}}{\left\|  -\dfrac A{\sqrt{2}}+\dfrac D{\sqrt{2}}\right\|
^{2}}>\dfrac13, &  & \dfrac{\left\|  -\sqrt{\dfrac23}A+\dfrac B{\sqrt{6}%
}+\dfrac D{\sqrt{6}}\right\|  ^{2}}{\left\|  -\dfrac B{\sqrt{2}}+\dfrac
D{\sqrt{2}}\right\|  ^{2}}>\dfrac13,\\
&  & \\
\dfrac{\left\|  \dfrac A{\sqrt{6}}-\sqrt{\dfrac23}C+\dfrac D{\sqrt{6}%
}\right\|  ^{2}}{\left\|  -\dfrac A{\sqrt{2}}+\dfrac D{\sqrt{2}}\right\|
^{2}}>\dfrac13, &  & \dfrac{\left\|  -\sqrt{\dfrac23}A+\dfrac C{\sqrt{6}%
}+\dfrac D{\sqrt{6}}\right\|  ^{2}}{\left\|  -\dfrac C{\sqrt{2}}+\dfrac
D{\sqrt{2}}\right\|  ^{2}}>\dfrac13.
\end{array}
\]
A\ computation of the probability that all four inequalities hold
simultaneously does not seem possible. We note that the two expressions in the
left-hand column possess the same denominator, which is essential for the next
section. An entirely different approach will be needed for the two expressions
in the right-hand column. Finally, the two left-hand inequalities are true if
and only if the orthogonal projection of $\vec w$ is $r\,\vec u+s\,\vec v$ for
some $r\leq1$, $s\leq1$; this translates into $\tilde D$ falling in the cone
\[
A+\left\{  (1-r)\,(B-A)+(1-s)\,(C-A):r\geq0,\,s\geq0\right\}  =(B+C)-\Gamma.
\]
Likewise, the two right-hand inequalities are true if and only if the
orthogonal projection of $\vec w$ is $r\,\vec u+s\,\vec v$ for some $r\geq0$,
$s\geq0$; this translates immediately into $\tilde D$ falling in the cone
$\Gamma$.

\section{Krishnaiah Bivariate $F$-Ratio Distribution}

Let $X$ be an $n\times2$ matrix of $n$ independent random row $2$-vectors,
each distributed according to $N(0,\Sigma)$, where
\[%
\begin{array}
[c]{ccc}%
\Sigma=\left(
\begin{array}
[c]{cc}%
1 & \rho\\
\rho & 1
\end{array}
\right)  , &  & -1<\rho<1.
\end{array}
\]
Let $\hat\sigma_{ij}=X_{i}\cdot X_{j}$, where $X_{k}$ is the $k^{\text{th}}$
column of $X$ for $k=1,2$. Also let $\hat\tau=Y\cdot Y$, where $Y$ is a random
column $m$-vector satisfying $Y\sim N(0,I)$ and independent of $X$. Then the
joint distribution of $(m\,\hat\sigma_{11})/(n\,\hat\tau)$, $(m\,\hat
\sigma_{22})/(n\,\hat\tau)$ is a bivariate $F$-ratio distribution with $(n,m)$
degrees of freedom and with $\Sigma$ as the associated covariance matrix (for
$X$). These conditions are clearly met in our case, for which $n=m=3$,
\[%
\begin{array}
[c]{ccc}%
X=\left(
\begin{array}
[c]{ccc}%
\dfrac{a_{1}}{\sqrt{6}}-\sqrt{\dfrac23}b_{1}+\dfrac{d_{1}}{\sqrt{6}} &  &
\dfrac{a_{1}}{\sqrt{6}}-\sqrt{\dfrac23}c_{1}+\dfrac{d_{1}}{\sqrt{6}}\\
\dfrac{a_{2}}{\sqrt{6}}-\sqrt{\dfrac23}b_{2}+\dfrac{d_{2}}{\sqrt{6}} &  &
\dfrac{a_{2}}{\sqrt{6}}-\sqrt{\dfrac23}c_{2}+\dfrac{d_{2}}{\sqrt{6}}\\
\dfrac{a_{3}}{\sqrt{6}}-\sqrt{\dfrac23}b_{3}+\dfrac{d_{3}}{\sqrt{6}} &  &
\dfrac{a_{3}}{\sqrt{6}}-\sqrt{\dfrac23}c_{3}+\dfrac{d_{3}}{\sqrt{6}}%
\end{array}
\right)  , &  & Y=\left(
\begin{array}
[c]{c}%
-\dfrac{a_{1}}{\sqrt{2}}+\dfrac{d_{1}}{\sqrt{2}}\\
-\dfrac{a_{2}}{\sqrt{2}}+\dfrac{d_{2}}{\sqrt{2}}\\
-\dfrac{a_{3}}{\sqrt{2}}+\dfrac{d_{3}}{\sqrt{2}}%
\end{array}
\right)
\end{array}
\]
and
\begin{align*}
\rho &  =\operatorname*{Cov}\left(  \dfrac{a_{k}}{\sqrt{6}}-\sqrt{\dfrac
23}b_{k}+\dfrac{d_{k}}{\sqrt{6}},\dfrac{a_{k}}{\sqrt{6}}-\sqrt{\dfrac23}%
c_{k}+\dfrac{d_{k}}{\sqrt{6}}\right) \\
&  =\operatorname*{Var}\left(  \dfrac{a_{k}}{\sqrt{6}}+\dfrac{d_{k}}{\sqrt{6}%
}\right)  =\frac13.
\end{align*}
It follows that \cite{KA}
\[
\operatorname*{P}\left(  \frac{m\,\hat\sigma_{11}}{n\,\hat\tau}>\xi\text{ and
}\frac{m\,\hat\sigma_{22}}{n\,\hat\tau}>\xi\right)  =\frac{(1-\rho^{2})^{n/2}%
}{\Gamma(m/2)\Gamma(n/2)}%
{\displaystyle\sum\limits_{k=0}^{\infty}}
\frac{\rho^{2k}\Gamma(n+m/2+2k)}{k!\,\Gamma(n/2+k)}\Lambda_{k}
\]
where
\[
\Lambda_{k}=%
{\displaystyle\int\limits_{\eta}^{\infty}}
{\displaystyle\int\limits_{\eta}^{\infty}}
\frac{(x\,y)^{n/2+k-1}}{(1+x+y)^{n+2k+m/2}}dx\,dy.
\]
and $\eta=\left(  n\,\xi\right)  /\left(  m(1-\rho^{2})\right)  $. Each
$\Lambda_{k}$ can be evaluated symbolically and the series appears to converge
fairly quickly. For our case, $\xi=1/3$, hence $\eta=3/8$ and the desired
probability is $0.6810669069...$.

\section{Miller Bivariate Density for Products}

We want the probability that both $(B-A)\cdot(D-A)>0$ and $(C-A)\cdot(D-A)>0.$
For now, we consider a simpler scenario in which $a=A,b=B,c=C,d=D$ are scalars.

Let $X$ be a random $p$-vector and $y$ be a random scalar. Assume that
$(X,y)\sim N(0,\Sigma)$ and that the $(p+1)\times(p+1)$ covariance matrix
$\Sigma$ has inverse
\[
\Sigma^{-1}=\left(
\begin{array}
[c]{cc}%
\Omega & \upsilon\\
\upsilon^{\prime} & \omega
\end{array}
\right)
\]
where $\Omega$ is $p\times p$, $\upsilon$ is $p\times1$ and $\omega$ is a
scalar. Let $Z=y\,X$. Then the joint density of $Z$ is \cite{Mi}
\[
\frac{2\sqrt{\det(\Sigma^{-1})}}{(2\pi)^{(p+1)/2}}\left(  \frac\omega
{Z^{\prime}\Omega Z}\right)  ^{(p-1)/4}\exp(-\upsilon^{\prime}Z)K_{(p-1)/2}%
\left(  \sqrt{\omega\,Z^{\prime}\Omega Z}\right)
\]
where $K_{(p-1)/2}(\theta)$ is the modified Bessel function of the second
kind. In our case, $p=2$,
\[
\Sigma=\operatorname*{Cov}\left(  \left(
\begin{array}
[c]{c}%
x_{1}\\
x_{2}\\
y
\end{array}
\right)  \right)  =\operatorname*{Cov}\left(  \left(
\begin{array}
[c]{c}%
b-a\\
c-a\\
d-a
\end{array}
\right)  \right)  =\left(
\begin{array}
[c]{ccc}%
2 & 1 & 1\\
1 & 2 & 1\\
1 & 1 & 2
\end{array}
\right)
\]
hence
\[%
\begin{array}
[c]{ccccccc}%
\Sigma^{-1}=\left(
\begin{array}
[c]{ccc}%
3/4 & -1/4 & -1/4\\
-1/4 & 3/4 & -1/4\\
-1/4 & -1/4 & 3/4
\end{array}
\right)  , &  & \Omega=\dfrac14\left(
\begin{array}
[c]{cc}%
3 & -1\\
-1 & 3
\end{array}
\right)  , &  & \upsilon=-\dfrac14\left(
\begin{array}
[c]{c}%
1\\
1
\end{array}
\right)  , &  & \omega=\dfrac34,
\end{array}
\]
\[%
\begin{array}
[c]{ccccc}%
\sqrt{\det(\Sigma^{-1})}=\dfrac12, &  & -\upsilon^{\prime}Z=\dfrac14\left(
z_{1}+z_{2}\right)  , &  & Z^{\prime}\Omega Z=\dfrac14\left(  3z_{1}%
^{2}-2z_{1}z_{2}+3z_{2}^{2}\right)  .
\end{array}
\]
Also
\[
K_{1/2}(\theta)=\sqrt{\frac\pi2}\frac{\exp(-\theta)}{\sqrt{\theta}}
\]
thus the density simplifies to
\[
f(z_{1},z_{2})=\frac1{2\pi}\frac{\exp\left(  \dfrac14\left(  z_{1}+z_{2}%
-\sqrt{3}\sqrt{3z_{1}^{2}-2z_{1}z_{2}+3z_{2}^{2}}\right)  \right)  }%
{\sqrt{3z_{1}^{2}-2z_{1}z_{2}+3z_{2}^{2}}}.
\]

In fact, we wish to compute the probability that the sum of \textit{three}
independent copies of $Z$ has both components $>0$. One way to do this is to
evaluate the sextuple integral:
\[%
{\displaystyle\int\limits_{0}^{\infty}}
{\displaystyle\int\limits_{0}^{\infty}}
{\displaystyle\int\limits_{-\infty}^{\infty}}
{\displaystyle\int\limits_{-\infty}^{\infty}}
{\displaystyle\int\limits_{-\infty}^{\infty}}
{\displaystyle\int\limits_{-\infty}^{\infty}}
f(z_{1}-u_{1}-v_{1},z_{2}-u_{2}-v_{2})f(u_{1},u_{2})f(v_{1},v_{2}%
)\,du_{1}\,du_{2}\,dv_{1}\,dv_{2}\,dz_{1}\,dz_{2}=0.683...
\]
obtained via repeated convolution of $f$ with itself. Another way uses the
Fourier transform
\[
F(w_{1},w_{2})=\frac1{\sqrt{(w_{1}-i)(3w_{1}+i)+(w_{2}-i)(3w_{2}%
+i)+2w_{1}w_{2}-1}}
\]
which, when cubed, maps back to a remarkably simple density function. Higher
precision is now possible:
\[%
{\displaystyle\int\limits_{0}^{\infty}}
{\displaystyle\int\limits_{0}^{\infty}}
\frac1{4\sqrt{3}\pi}\exp\left(  \dfrac14\left(  z_{1}+z_{2}-\sqrt{3}%
\sqrt{3z_{1}^{2}-2z_{1}z_{2}+3z_{2}^{2}}\right)  \right)  dz_{1}%
\,dz_{2}=0.6837762984...
\]
and the elaborate details appear later.

It is surprising that working with $(B-A)\cdot(D-A)$, $(C-A)\cdot(D-A)$
simultaneously should be so difficult. The two expressions are familiar: they
are sample covariance coefficients $\hat\gamma_{13}$, $\hat\gamma_{23}$
respectively between samples of size $=3$. A marginal density for either is
found in \cite{Pr1}, but not much else is known about off-diagonal elements of
a Wishart matrix \cite{Pr2}. Of course, $\hat\gamma_{13}>0$ and $\hat
\gamma_{23}>0$ if and only if the corresponding sample correlation
coefficients $\hat\rho_{13}>0$ and $\hat\rho_{23}>0$. A formula for
a\ trivariate density for $(\hat\rho_{12},\hat\rho_{13},\hat\rho_{23})$ is
outlined in \cite{Fs, Be} -- evidently a sample size $>4$ is presumed -- and
details still need to come together.

\section{Pinned Simplices}

A slight variation on defining a random tetrahedron in $3$-space is to keep
one vertex fixed at the origin and to select the other three vertices
independently from $N(0,I)$ as before. We say that the tetrahedron is
\textbf{pinned}.

For the moment, let us consider pinned random Gaussian triangles $ABC$ in
$2$-space with $C=(0,0)$. Note that $(B-A)\cdot(-A)$ is equal to \cite{EiS}
\begin{align*}
&  \ \ \frac{1+\sqrt{2}}2\left\|  -\frac{\sqrt{2+\sqrt{2}}}2A+\frac
{\sqrt{2-\sqrt{2}}}2B\right\|  ^{2}-\ \frac{-1+\sqrt{2}}2\left\|  \frac
{\sqrt{2-\sqrt{2}}}2A+\frac{\sqrt{2+\sqrt{2}}}2B\right\|  ^{2}\\
\  &  =\ \frac{1+\sqrt{2}}2\chi_{2}^{2}-\frac{-1+\sqrt{2}}2\bar\chi_{2}^{2},
\end{align*}
therefore
\[
p=\operatorname*{P}\left(  \frac{1+\sqrt{2}}2\chi_{2}^{2}-\frac{-1+\sqrt{2}%
}2\bar\chi_{2}^{2}>0\right)  =\operatorname*{P}\left(  F_{2,2}>3-2\sqrt
{2}\right)  =\frac{2+\sqrt{2}}4.
\]
It follows that $\tilde C$ falls between $A$ and $B$ on $Q$ with probability
$2p-1=1/\sqrt{2}$. Equivalently, $\operatorname*{P}(\alpha>\pi/2)=1-p$ because
$\cos(\alpha)<0$ if and only if $(B-A)\cdot(-A)<0$. By symmetry,
$\operatorname*{P}(\beta>\pi/2)=1-p$ as well and $\operatorname*{P}(\gamma
>\pi/2)=1/2$. We deduce that
\begin{align*}
&  \operatorname*{P}\left(  \text{pinned triangle }ABC\text{ is acute}\right)
\\
&  =1-\left(  \operatorname*{P}(\alpha>\pi/2)+\operatorname*{P}(\beta
>\pi/2)+\operatorname*{P}(\gamma>\pi/2)\right) \\
\  &  =1-\left(  2-2p+1/2\right)  =-1/2+1/\sqrt{2}%
\end{align*}
as was to be shown.

Let us now return to pinned random Gaussian tetrahedra in $3$-space with
$D=(0,0,0)$. For brevity, we focus only on the probability that both
$(B-A)\cdot(-A)>0$ and $(C-A)\cdot(-A)>0.$ Consider a simpler scenario in
which $a=A,b=B,c=C$ are scalars. Using Miller's \cite{Mi} formulas,
\[
\Sigma=\operatorname*{Cov}\left(  \left(
\begin{array}
[c]{c}%
b-a\\
c-a\\
-a
\end{array}
\right)  \right)  =\left(
\begin{array}
[c]{ccc}%
2 & 1 & 1\\
1 & 2 & 1\\
1 & 1 & 1
\end{array}
\right)
\]
hence
\[%
\begin{array}
[c]{ccccccc}%
\Sigma^{-1}=\left(
\begin{array}
[c]{ccc}%
1 & 0 & -1\\
0 & 1 & -1\\
-1 & -1 & 3
\end{array}
\right)  , &  & \Omega=\left(
\begin{array}
[c]{cc}%
1 & 0\\
0 & 1
\end{array}
\right)  , &  & \upsilon=-\left(
\begin{array}
[c]{c}%
1\\
1
\end{array}
\right)  , &  & \omega=3,
\end{array}
\]
\[%
\begin{array}
[c]{ccccc}%
\sqrt{\det(\Sigma^{-1})}=1, &  & -\upsilon^{\prime}Z=z_{1}+z_{2}, &  &
Z^{\prime}\Omega Z=z_{1}^{2}+z_{2}^{2},
\end{array}
\]
\[
f(z_{1},z_{2})=\frac1{2\pi}\frac{\exp\left(  z_{1}+z_{2}-\sqrt{3}\sqrt
{z_{1}^{2}+z_{2}^{2}}\right)  }{\sqrt{z_{1}^{2}+z_{2}^{2}}}.
\]
The probability we want is given by the sextuple integral, which has value
$0.834...$, but can be computed more accurately via the Fourier transform
\[
F(w_{1},w_{2})=\frac1{\sqrt{(w_{1}-i)^{2}+(w_{2}-i)^{2}+3}}.
\]
Mapping $F(w_{1},w_{2})^{3}$ from frequency back to signal domain, we
calculate
\[%
{\displaystyle\int\limits_{0}^{\infty}}
{\displaystyle\int\limits_{0}^{\infty}}
\frac1{2\sqrt{3}\pi}\exp\left(  z_{1}+z_{2}-\sqrt{3}\sqrt{z_{1}^{2}+z_{2}^{2}%
}\right)  dz_{1}\,dz_{2}=0.8343764256....
\]
After a discussion of some related problems in geometric probability
\cite{F1}, more information on $F(w_{1},w_{2})$ (also called the
\textit{characteristic function} corresponding to $f(z_{1},z_{2})$) will be given.

\section{Random Solid Angles}

We restrict attention to pinned Gaussian random tetrahedra $ABCD$ with
$D=(0,0,0)$. The \textbf{dihedral angle} $\alpha$ is the angle between normal
vectors $A\times B$, $A\times C$ to the triangular faces $ADB$, $ADC$
respectively:
\[
\alpha=\arccos\left(  \frac{(A\times B)\cdot(A\times C)}{\left\|  A\times
B\right\|  \,\left\|  A\times C\right\|  }\right)  .
\]
Angles $\beta$ and $\gamma$ are defined likewise. For example, a
\textit{regular} tetrahedron has dihedral angles each equal to
\[
\arccos(1/3)=1.2309594173...=\pi-1.9106332362...\approx70.53^{\circ}%
\approx180^{\circ}-109.47^{\circ}.
\]
The joint density for $\alpha$, $\beta$, $\gamma$ is \cite{Ms}
\[
\left\{
\begin{array}
[c]{l}%
-\dfrac1\pi\dfrac{\cos\left(  \dfrac{x+y+z}2\right)  \cos\left(
\dfrac{-x+y+z}2\right)  \cos\left(  \dfrac{x-y+z}2\right)  \cos\left(
\dfrac{x+y-z}2\right)  }{\sin(x)^{2}\sin(y)^{2}\sin(z)^{2}}\\
\;\;\;\;\;\;\;\;\;\;\;\;\;\;\;\text{if }x+y+z>\pi\text{, }x+y<\pi+z\text{,
}y+z<\pi+x\text{ and }z+x<\pi+y,\\
0\;\;\;\;\;\;\;\;\;\;\;\;\;\text{otherwise.}%
\end{array}
\right.
\]
As a consequence, $\alpha$ is uniformly distributed on $[0,\pi]$ and angles
$\alpha$, $\beta$, $\gamma$ are uncorrelated but pairwise dependent
($(\alpha,\beta)$ is not uniform on $[0,\pi]\times[0,\pi]$).

The \textbf{solid angle} (or \textbf{trihedral angle}) at $D$ is the area of
the spherical triangle on the unit sphere, center $D$, with vertices
$A/\left\|  A\right\|  $, $B/\left\|  B\right\|  $, $C/\left\|  C\right\|  $.
It is equal to the spherical excess $\alpha+\beta+\gamma-\pi$, which is
between $0$ and $2\pi$. It is also equal to \cite{OS, Er}
\[
\left\{
\begin{array}
[c]{lll}%
2\arctan(\zeta) &  & \text{if }\zeta\geq0,\\
2\pi+2\arctan(\zeta) &  & \text{if }\zeta<0
\end{array}
\right.
\]
where
\[
\zeta=\frac{\left|  A\cdot(B\times C)\right|  }{\left\|  A\right\|  \,\left\|
B\right\|  \,\left\|  C\right\|  +(A\cdot B)\left\|  C\right\|  +(A\cdot
C)\left\|  B\right\|  +(B\cdot C)\left\|  A\right\|  }.
\]
It can be shown that, if a tetrahedron is acute, then each of its four solid
angles are less than $\pi/2$, but not conversely \cite{KS}. A\ proposed
density for the solid angle at $D$ was published in 1867 \cite{CE}:
\[
-\frac{(x^{2}-4\pi x+3\pi^{2}-6)\cos(x)-6(x-2\pi)\sin(x)-2(x^{2}-4\pi
x+3\pi^{2}+3)}{16\pi\cos(x/2)^{4}}
\]
for $0<x<2\pi$ and remained obscure until it was cited in a recent paper
\cite{HDB}. Details of the supporting geometric proof need to be carefully
examined. No analytic proof using the joint density for $\alpha,\beta,\gamma$
has yet been found.

As far as is known, no analogous results are known for general Gaussian random
tetrahedra. In particular, the sum $\sigma$ of the four solid angles
associated with a tetrahedra $T$ possesses a fascinating property \cite{FK}:
\[
\frac\sigma{2\pi}=\operatorname*{P}\left(
\begin{array}
[c]{c}%
\text{the orthogonal projection of }T\text{ onto a uniform}\\
\text{random plane in }3\text{-space is a triangle}%
\end{array}
\right)
\]
and it would be good to understand $\sigma$ more fully. As an example, the
regular tetrahedra has solid angles each equal to
\[%
\begin{array}
[c]{ccc}%
3\arccos(1/3)-\pi=0.5512855984..., &  & \text{hence }\sigma/(2\pi
)=0.3509593121...
\end{array}
\]
and this is the maximum such value over all \textit{equifacial} tetrahedra
(all faces are congruent) \cite{Ka, Hv, KM, Hj}. No one has studied the
distribution of $\sigma$ when $T$ is itself allowed to be random.

While we know the mean volume of a tetrahedron \cite{Hs, Sm, Zn, Ph} with
uniform random vertices in the unit ball ($12\pi/715$) and with uniform random
vertices in the unit cube ($3977/21600-\pi^{2}/2160$), the Gaussian random
scenario remains open (there is doubt about claims in \cite{Bo}).

In $2$-space, a triangle is acute if and only if its circumcenter lies inside
the triangle (the circumcircle contains all three vertices). In $3$-space, a
tetrahedron is $3$\textbf{-well-centered} if its circumcenter lies inside the
tetrahedron; a tetrahedron is $2$\textbf{-well-centered} if the circumcenter
of each face lies inside the face. An acute $T$ can fail to be $3$%
-well-centered, and a $3$-well-centered $T$ can fail to be acute. However, an
acute $T$ must be $2$-well-centered (equivalently, all its faces must be
acute) but not conversely \cite{ESU, VZH}. Many problems involving Gaussian
random tetrahedra suggest themselves.

\section{Fourier Transforms}

It is not difficult to prove all the formulas we need in the ``forward''
direction (starting with $f$ and ending with $F$). This is done first for the
pinned case, which is easier, and then for the general case. Motivating the
formula for the inverse Fourier transform of $F^{3}$ is harder. We do this for
the pinned case only.

\subsection{Pinned Case:\ Forward Direction}

Our objective is to evaluate two integrals:
\[
F(u,v)=\dfrac1{2\pi}%
{\displaystyle\int\limits_{-\infty}^{\infty}}
{\displaystyle\int\limits_{-\infty}^{\infty}}
\dfrac{\exp\left(  (1+i\,u)x+(1+i\,v)y-\sqrt{3}\sqrt{x^{2}+y^{2}}\right)
}{\sqrt{x^{2}+y^{2}}}dx\,dy,
\]
\[
G(u,v)=\dfrac1{2\sqrt{3}\pi}%
{\displaystyle\int\limits_{-\infty}^{\infty}}
{\displaystyle\int\limits_{-\infty}^{\infty}}
\exp\left(  (1+i\,u)x+(1+i\,v)y-\sqrt{3}\sqrt{x^{2}+y^{2}}\right)  dx\,dy.
\]
Let $x=r\cos(\theta)$, $y=r\sin(\theta)$, then $dx\,dy=r\,dr\,d\theta$ and
\begin{align*}
F(u,v)  &  =\dfrac1{2\pi}%
{\displaystyle\int\limits_{0}^{2\pi}}
{\displaystyle\int\limits_{0}^{\infty}}
\exp\left(  r\left[  (1+i\,u)\cos(\theta)+(1+i\,v)\sin(\theta)-\sqrt
{3}\right]  \right)  dr\,d\theta\\
\  &  =-\dfrac1{2\pi}%
{\displaystyle\int\limits_{0}^{2\pi}}
\frac1{(1+i\,u)\cos(\theta)+(1+i\,v)\sin(\theta)-\sqrt{3}}d\theta,
\end{align*}
\begin{align*}
G(u,v)  &  =\dfrac1{2\sqrt{3}\pi}%
{\displaystyle\int\limits_{0}^{2\pi}}
{\displaystyle\int\limits_{0}^{\infty}}
r\exp\left(  r\left[  (1+i\,u)\cos(\theta)+(1+i\,v)\sin(\theta)-\sqrt
{3}\right]  \right)  dr\,d\theta\\
\  &  =\dfrac1{2\sqrt{3}\pi}%
{\displaystyle\int\limits_{0}^{2\pi}}
\frac1{\left[  (1+i\,u)\cos(\theta)+(1+i\,v)\sin(\theta)-\sqrt{3}\right]
^{2}}d\theta
\end{align*}
since
\begin{align*}
\operatorname*{Re}\left[  (1+i\,u)\cos(\theta)+(1+i\,v)\sin(\theta)-\sqrt
{3}\right]   &  =\cos(\theta)+\sin(\theta)-\sqrt{3}\\
\  &  \leq\sqrt{2}-\sqrt{3}<0.
\end{align*}
Let $z=\exp(i\,\theta)$, then $d\theta=-i\,dz/z$ and
\begin{align*}
F(u,v)  &  =\dfrac i{2\pi}%
{\displaystyle\int\nolimits_{C}}
\frac1{(1+i\,u)\frac12(z+\frac1z)+(1+i\,v)\frac1{2i}(z-\frac1z)-\sqrt{3}}%
\frac{dz}z\\
&  =\dfrac i\pi%
{\displaystyle\int\nolimits_{C}}
\frac1{(1+i\,u)(z^{2}+1)-i(1+i\,v)(z^{2}-1)-2\sqrt{3}z}dz\\
&  =\dfrac i\pi%
{\displaystyle\int\nolimits_{C}}
\frac1{(v+1+i\,u-i)z^{2}-2\sqrt{3}z-(v-1-i\,u-i)}dz,
\end{align*}
\begin{align*}
G(u,v)  &  =-\dfrac i{2\sqrt{3}\pi}%
{\displaystyle\int\nolimits_{C}}
\frac1{\left[  (1+i\,u)\frac12(z+\frac1z)+(1+i\,v)\frac1{2i}(z-\frac
1z)-\sqrt{3}\right]  ^{2}}\frac{dz}z\\
&  =-\dfrac{2i}{\sqrt{3}\pi}%
{\displaystyle\int\nolimits_{C}}
\frac z{\left[  (v+1+i\,u-i)z^{2}-2\sqrt{3}z-(v-1-i\,u-i)\right]  ^{2}}dz
\end{align*}
where $C$ denotes the unit circle, center $0$, in the complex plane. The two
poles $z_{\text{pos}}$, $z_{\text{neg}}$ of each integrand are
\[
\frac{2\sqrt{3}\pm\sqrt{12+4(v+1+i\,u-i)(v-1-i\,u-i)}}{2(v+1+i\,u-i)}%
=\frac{\sqrt{3}\pm\sqrt{(u-i)^{2}+(v-i)^{2}+3}}{v+1+i\,u-i}
\]
and $z_{\text{neg}}$ is always inside $C$, $z_{\text{pos}}$ is always outside.
For $F$, $z_{\text{neg}}$ is a pole of order $1$ and the associated residue
is
\[
\lim_{z\rightarrow z_{\text{neg}}}\frac1{(v+1+i\,u-i)(z-z_{\text{pos}}%
)}=-\frac12\frac1{\sqrt{(u-i)^{2}+(v-i)^{2}+3}};
\]
multiplying by $(2\pi i)(i/\pi)$ completes the proof. For $G$, $z_{\text{neg}%
}$ is a pole of order $2$ and the associated residue is
\[
\lim_{z\rightarrow z_{\text{neg}}}\frac d{dz}\left\{  \frac z{(v+1+i\,u-i)^{2}%
(z-z_{\text{pos}})^{2}}\right\}  =\frac{\sqrt{3}}4\frac1{\left[
(u-i)^{2}+(v-i)^{2}+3\right]  ^{3/2}};
\]
multiplying by $(2\pi i)(-2i/(\sqrt{3}\pi))$ completes the proof.

\subsection{General Case:\ Forward Direction}

Our objective is to evaluate two integrals:
\[
F(u,v)=\dfrac{1}{2\pi}%
{\displaystyle\int\limits_{-\infty}^{\infty}}
{\displaystyle\int\limits_{-\infty}^{\infty}}
\dfrac{\exp\left(  \dfrac{1}{4}\left(  (1+4i\,u)x+(1+4i\,v)y-\sqrt{3}%
\sqrt{3x^{2}-2x\,y+3y^{2}}\right)  \right)  }{\sqrt{3x^{2}-2x\,y+3y^{2}}%
}dx\,dy,
\]%
\[
G(u,v)=\dfrac{1}{4\sqrt{3}\pi}%
{\displaystyle\int\limits_{-\infty}^{\infty}}
{\displaystyle\int\limits_{-\infty}^{\infty}}
\exp\left(  \dfrac{1}{4}\left(  (1+4i\,u)x+(1+4i\,v)y-\sqrt{3}\sqrt
{3x^{2}-2x\,y+3y^{2}}\right)  \right)  dx\,dy.
\]
Let
\[%
\begin{array}
[c]{ccc}%
x=\dfrac{r}{2\sqrt{2}}\left(  \sqrt{2}\cos(\theta)-\sin(\theta)\right)  , &  &
y=\dfrac{r}{2\sqrt{2}}\left(  \sqrt{2}\cos(\theta)+\sin(\theta)\right)
\end{array}
\]
then $3x^{2}-2x\,y+3y^{2}=r^{2}$ and the Jacobian determinant is
\[
\left(  \dfrac{1}{2\sqrt{2}}\right)  ^{2}\left\vert
\begin{array}
[c]{ccc}%
\sqrt{2}\cos(\theta)-\sin(\theta) &  & r\left(  -\sqrt{2}\sin(\theta
)-\cos(\theta)\right) \\
\sqrt{2}\cos(\theta)+\sin(\theta) &  & r\left(  -\sqrt{2}\sin(\theta
)+\cos(\theta)\right)
\end{array}
\right\vert =\frac{r}{2\sqrt{2}}%
\]
hence $dx\,dy=(r/(2\sqrt{2})\,dr\,d\theta$. We obtain
\begin{align*}
F(u,v)  &  =\dfrac{1}{2\pi}\dfrac{1}{2\sqrt{2}}%
{\displaystyle\int\limits_{0}^{2\pi}}
{\displaystyle\int\limits_{0}^{\infty}}
\exp\left(  \frac{1}{4}\dfrac{r}{2\sqrt{2}}\left[  \left(  1+4i\,u\right)
\left(  \sqrt{2}\cos(\theta)-\sin(\theta)\right)  \right.  \right. \\
&  \;\;\;\;\;\;\;\;\;\;\;\;\;\left.  \left.  +\left(  1+4i\,v\right)  \left(
\sqrt{2}\cos(\theta)+\sin(\theta)\right)  -2\sqrt{2}\cdot\sqrt{3}\right]
\right)  dr\,d\theta\\
&  =\frac{1}{4\sqrt{2}\pi}%
{\displaystyle\int\limits_{0}^{2\pi}}
{\displaystyle\int\limits_{0}^{\infty}}
\exp\left(  \frac{r}{4}\left[  \cos(\theta)+2i(u+v)\cos(\theta)-\sqrt
{2}i(u-v)\sin(\theta)-\sqrt{3}\right]  \right)  dr\,d\theta\\
\  &  =-\dfrac{1}{\sqrt{2}\pi}%
{\displaystyle\int\limits_{0}^{2\pi}}
\frac{1}{\cos(\theta)+2i(u+v)\cos(\theta)-\sqrt{2}i(u-v)\sin(\theta)-\sqrt{3}%
}d\theta,
\end{align*}%
\begin{align*}
G(u,v)  &  =\dfrac{1}{4\sqrt{3}\pi}\dfrac{1}{2\sqrt{2}}%
{\displaystyle\int\limits_{0}^{2\pi}}
{\displaystyle\int\limits_{0}^{\infty}}
r\exp\left(  \frac{1}{4}\dfrac{r}{2\sqrt{2}}\left[  \left(  1+4i\,u\right)
\left(  \sqrt{2}\cos(\theta)-\sin(\theta)\right)  \right.  \right. \\
&  \ \;\;\;\;\;\;\;\;\;\;\;\;\;\left.  \left.  +\left(  1+4i\,v\right)
\left(  \sqrt{2}\cos(\theta)+\sin(\theta)\right)  -2\sqrt{2}\cdot\sqrt
{3}\right]  \right)  dr\,d\theta\\
\  &  =\frac{1}{8\sqrt{6}\pi}%
{\displaystyle\int\limits_{0}^{2\pi}}
{\displaystyle\int\limits_{0}^{\infty}}
r\exp\left(  \frac{r}{4}\left[  \cos(\theta)+2i(u+v)\cos(\theta)-\sqrt
{2}i(u-v)\sin(\theta)-\sqrt{3}\right]  \right)  dr\,d\theta\\
\  &  =\dfrac{2}{\sqrt{6}\pi}%
{\displaystyle\int\limits_{0}^{2\pi}}
\frac{1}{\left[  \cos(\theta)+2i(u+v)\cos(\theta)-\sqrt{2}i(u-v)\sin
(\theta)-\sqrt{3}\right]  ^{2}}d\theta
\end{align*}
since
\begin{align*}
\operatorname*{Re}\left[  \cos(\theta)+2i(u+v)\cos(\theta)-\sqrt{2}%
i(u-v)\sin(\theta)-\sqrt{3}\right]   &  =\cos(\theta)-\sqrt{3}\\
&  \leq1-\sqrt{3}<0.
\end{align*}
Let $z=\exp(i\,\theta)$, then $d\theta=-i\,dz/z$ and
\begin{align*}
F(u,v)  &  =\dfrac{i}{\sqrt{2}\pi}%
{\displaystyle\int\nolimits_{C}}
\frac{1}{\frac{1}{2}(z+\frac{1}{z})+2i(u+v)\frac{1}{2}(z+\frac{1}{z})-\sqrt
{2}i(u-v)\frac{1}{2i}(z-\frac{1}{z})-\sqrt{3}}\frac{dz}{z}\\
\  &  =\dfrac{\sqrt{2}i}{\pi}%
{\displaystyle\int\nolimits_{C}}
\frac{1}{(z^{2}+1)+2i(u+v)(z^{2}+1)-\sqrt{2}(u-v)(z^{2}-1)-2\sqrt{3}z}dz\\
\  &  =\dfrac{\sqrt{2}i}{\pi}%
{\displaystyle\int\nolimits_{C}}
\frac{1}{(1-\sqrt{2}u+\sqrt{2}v+2iu+2iv)z^{2}-2\sqrt{3}z+(1+\sqrt{2}u-\sqrt
{2}v+2iu+2iv)}dz,
\end{align*}%
\begin{align*}
G(u,v)  &  =-\dfrac{2i}{\sqrt{6}\pi}%
{\displaystyle\int\nolimits_{C}}
\frac{1}{\left[  \frac{1}{2}(z+\frac{1}{z})+2i(u+v)\frac{1}{2}(z+\frac{1}%
{z})-\sqrt{2}i(u-v)\frac{1}{2i}(z-\frac{1}{z})-\sqrt{3}\right]  ^{2}}\frac
{dz}{z}\\
&  =-\dfrac{8i}{\sqrt{6}\pi}%
{\displaystyle\int\nolimits_{C}}
\frac{z}{\left[  (z^{2}+1)+2i(u+v)(z^{2}+1)-\sqrt{2}(u-v)(z^{2}-1)-2\sqrt
{3}z\right]  ^{2}}dz\\
&  =-\dfrac{8i}{\sqrt{6}\pi}%
{\displaystyle\int\nolimits_{C}}
\frac{z}{\left[  (1-\sqrt{2}u+\sqrt{2}v+2iu+2iv)z^{2}-2\sqrt{3}z+(1+\sqrt
{2}u-\sqrt{2}v+2iu+2iv)\right]  ^{2}}dz.
\end{align*}
The two poles $z_{\text{pos}}$, $z_{\text{neg}}$ of each integrand are
\begin{align*}
&  \frac{2\sqrt{3}\pm\sqrt{12-4(1-\sqrt{2}u+\sqrt{2}v+2iu+2iv)(1+\sqrt
{2}u-\sqrt{2}v+2iu+2iv)}}{2(1-\sqrt{2}u+\sqrt{2}v+2iu+2iv)}\\
&  =\frac{\sqrt{3}\pm\sqrt{2}\sqrt{(u-i)(3u+i)+(v-i)(3v+i)+2uv-1}}{1-\sqrt
{2}u+\sqrt{2}v+2iu+2iv}%
\end{align*}
and $z_{\text{neg}}$ is always inside $C$, $z_{\text{pos}}$ is always outside.
For $F$, $z_{\text{neg}}$ is a pole of order $1$ and the associated residue
is
\begin{align*}
&  \lim_{z\rightarrow z_{\text{neg}}}\frac{1}{(1-\sqrt{2}u+\sqrt
{2}v+2iu+2iv)(z-z_{\text{pos}})}\\
&  =-\frac{1}{2\sqrt{2}}\frac{1}{\sqrt{(u-i)(3u+i)+(v-i)(3v+i)+2uv-1}};
\end{align*}
multiplying by $(2\pi i)(\sqrt{2}i/\pi)$ completes the proof. For $G$,
$z_{\text{neg}}$ is a pole of order $2$ and the associated residue is
\begin{align*}
&  \lim_{z\rightarrow z_{\text{neg}}}\frac{d}{dz}\left\{  \frac{z}{(1-\sqrt
{2}u+\sqrt{2}v+2iu+2iv)^{2}(z-z_{\text{pos}})^{2}}\right\} \\
&  =\frac{\sqrt{6}}{16}\frac{1}{\left[  (u-i)(3u+i)+(v-i)(3v+i)+2uv-1\right]
^{3/2}};
\end{align*}
multiplying by $(2\pi i)(-8i/(\sqrt{6}\pi))$ completes the proof. Clearly
there is common structure to both cases and a more encompassing theorem should
be possible.

\subsection{Pinned Case:\ Backward Direction}

The fact that $G=F^{3}$ is fairly miraculous but completely unmotivated. Let
us briefly sketch a \textquotedblleft backward\textquotedblright\ argument for
the pinned case only, starting with a known inverse Fourier transform
\cite{GR}:
\[%
{\displaystyle\int\limits_{-\infty}^{\infty}}
\frac{\exp(-iv\,y)}{\sqrt{(u-i)^{2}+(v-i)^{2}+3}}dv=2\exp(y)K_{0}\left(
|y|\sqrt{(u-i)^{2}+3}\right)  .
\]
Differentiating both sides with respect to $u$, we obtain
\[
-%
{\displaystyle\int\limits_{-\infty}^{\infty}}
\frac{\exp(-iv\,y)(u-i)}{\left[  (u-i)^{2}+(v-i)^{2}+3\right]  ^{3/2}%
}dv=-2\exp(y)K_{1}\left(  |y|\sqrt{(u-i)^{2}+3}\right)  \cdot\frac
{|y|(u-i)}{\sqrt{(u-i)^{2}+3}}
\]
that is,
\[%
{\displaystyle\int\limits_{-\infty}^{\infty}}
\frac{\exp(-iv\,y)}{\left[  (u-i)^{2}+(v-i)^{2}+3\right]  ^{3/2}}%
dv=2|y|\exp(y)\frac{K_{1}\left(  |y|\sqrt{(u-i)^{2}+3}\right)  }%
{\sqrt{(u-i)^{2}+3}}.
\]
Another known inverse Fourier transform is useful now \cite{GR}:
\begin{align*}
&
{\displaystyle\int\limits_{-\infty}^{\infty}}
\exp(-iu\,x)\frac{K_{1}\left(  |y|\sqrt{(u-i)^{2}+3}\right)  }{\sqrt
{(u-i)^{2}+3}}du\\
&  =\frac{\sqrt{2\pi}}{3^{1/4}}\frac{\exp(x)}{|y|}\left(  x^{2}+y^{2}\right)
^{1/4}K_{1/2}\left(  \sqrt{3}\sqrt{x^{2}+y^{2}}\right) \\
\  &  =\frac{\sqrt{2\pi}}{3^{1/4}}\frac{\exp(x)}{|y|}\left(  x^{2}%
+y^{2}\right)  ^{1/4}\sqrt{\frac{\pi}{2}}\frac{\exp\left(  -\sqrt{3}%
\sqrt{x^{2}+y^{2}}\right)  }{3^{1/4}\left(  x^{2}+y^{2}\right)  ^{1/4}}\\
\  &  =\frac{\pi}{\sqrt{3}}\frac{\exp(x)}{|y|}\exp\left(  -\sqrt{3}\sqrt
{x^{2}+y^{2}}\right)  .
\end{align*}
Multiplying both sides by $2|y|\exp(y)/(2\pi)^{2}$, we conclude that
\[
\frac{1}{(2\pi)^{2}}%
{\displaystyle\int\limits_{-\infty}^{\infty}}
{\displaystyle\int\limits_{-\infty}^{\infty}}
\frac{\exp(-iu\,x-iv\,y)}{\left[  (u-i)^{2}+(v-i)^{2}+3\right]  ^{3/2}%
}du\,dv=\frac{1}{2\sqrt{3}\pi}\exp\left(  x+y-\sqrt{3}\sqrt{x^{2}+y^{2}%
}\right)
\]
as was to be shown.

\section{Addendum}

The mean volume of an unpinned Gaussian random tetrahedron is $2\sqrt{2}%
/(3\pi)$. \ We do not know the corresponding result for pinned tetrahedra nor
any higher moments \cite{Efrn, MCRF}. \ See \cite{F2} for generalization of
Section 6 and \cite{F3} for proof of the Crofton-Exhumatus density.

\section{Acknowledgement}

I am grateful to Robert Israel for a helpful discussion about residue
calculus. Much more relevant material can be found at \cite{F4, F5}, including
experimental computer runs that aided theoretical discussion here.

\end{document}